\documentclass[12pt,
]{article}

\usepackage{amssymb}

\usepackage{amsmath}
\usepackage{mathrsfs}
\usepackage{theorem}
\usepackage{hyperref}

\usepackage{newtxtext}

\newcommand{\al}{\alpha}
\newcommand{\be}{\beta}
\newcommand{\ep}{\varepsilon}
\newcommand{\ga}{\gamma}
\newcommand{\la}{\lambda}
\newcommand{\om}{\omega}
\newcommand{\ro}{\varrho}
\newcommand{\de}{\partial}
\newcommand{\si}{\sigma}
\newcommand{\vf}{\varphi}
\newcommand{\De}{\Delta}
\newcommand\A{\mathop{\mathscr A}\nolimits}
\newcommand\B{\mathop{\mathscr B}\nolimits}
\newcommand{\Om}{\Omega}
\newcommand{\Si}{\Sigma}
\newcommand{\R}{\mathbb R}
\newcommand{\Z}{\mathbb Z}
\newcommand{\N}{\mathbb N}

\newcommand{\Ga}{\Gamma}

\newcommand\Cspt{\mathaccent"017{C}}
\newcommand{\dive}{\mathop{\rm div}\nolimits}
\newcommand{\curl}{\operatorname{curl}}
\newcommand{\w}{\wedge}

\renewcommand\leq{\leqslant}
\renewcommand\geq{\geqslant}

\newtheorem{theorem}{Theorem}
\newenvironment{proof} {{\bf Proof.}}{\hfill \fbox{}\\ \smallskip}
{\theorembodyfont{\rmfamily} \newtheorem{remark}{Remark}}

\title{New completeness theorems on the boundary in Elasticity}

\author{A. Cialdea
\thanks{Department of Mathematics, Computer Sciences and Economics,
University of Ba\-si\-li\-ca\-ta, V.le dell'Ateneo Lucano, 10, 85100 Potenza, Italy.
 \textit{email:}
alberto.cialdea@unibas.it. 
}
}
\date{}    

\begin{document}

{\footnotesize  © 2024. This manuscript version is made available under the CC-BY-NC-ND 4.0 license 
\url{https://creativecommons.org/licenses/by-nc-nd/4.0}
}

{\let\newpage\relax\maketitle}

\noindent\textbf{Abstract}
The completeness on the boundary (in the sense of Picone) of certain systems
related to the III and IV BVPs for the elasticity system is proved. The completeness is obtained in both $L^p$ ($1\leq 1<\infty$) and uniform norms.

 \section{Introduction}
 
 In the main boundary value problems in the classical theory of elasticity,  
 the displacements or the stresses on the boundary of the body are prescribed.
 Following the terminology 
 adopted in Kupradze's monograph \cite{kupra}, we call these problems the I and the II boundary value problem (BVP), respectively.
 
 Other interesting BVPs consider  the normal or the tangential components of the displacement or 
 the stress
 (see, e.g., Section VI in \cite{pragsyng47}).
 Here we consider the III and the IV BVP
 (according to \cite{kupra}).   In the third problem the normal component of the displacement vector and the tangential components of the stress are given on the boundary, while in the fourth problem, conversely, 
 the normal component of the stress and the tangential components of the displacement are given. 
 
 The present paper aims to study the completeness of certain systems on the boundary in different norms. 
 These  
 systems are constructed using polynomial solutions
 of the classical linear system of elasticity (which we call elastic polynomials) and are related to the III and IV BVPs.
  More specifically, let us denote by $T$ the stress operator on the boundary of a bounded domain
 and by $\{p_k\}$ the system of elastic polynomials. We determine the closure in $L^p$ and $C^0$ norms of the linear spaces
 spanned by $\left\{   p_k \cdot \nu \, , Tp_k - (Tp_k \cdot \nu)\nu
\right\}$ and $\left\{   p_k - (p_k \cdot \nu)\nu  \, , Tp_k \cdot \nu
\right\}$, $\nu$ being the  (outward) unit normal.
 
 This problem falls into the area of the so-called completeness theorems in the sense of Picone. 
 For an introduction to this topic and some comparison with completeness theorems 
of Mergelyan and Runge type we refer to \cite[pp.95--100]{C2019}.
Concerning the completeness in the sense of Picone of polynomial solutions, there are several available results for scalar partial differential equations. We refer to \cite{C2019} for a complete list of references on this subject. Here we mention that, concerning  the Dirichlet problem, necessary and sufficient conditions for
the completeness of polynomial solutions on the boundary in $L^p$ and $C^0$ norms are known 
in the case of scalar elliptic partial differential operators with constant coefficients \cite{C2007,C2012}.
As far as systems are concerned,
 very little is known. There are available results only for elasticity \cite{C19882,Fichera50}, thermoelasticity \cite{C2019} and Stokes system \cite{CiNino}.
 
 We remark that generally speaking, an existence and uniqueness
 theorem for a certain BVP and the existence of polynomial solutions of the partial differential equation do not necessarily imply the corresponding
 completeness on the boundary. An example is given by the Dirichlet problem for the operator $\De^2-\De$ (see \cite[p.100]{C2019}).

As Fichera writes in \cite{fichera1979} in the case of Laplace equation, the completeness in the sense of Picone is, \textit{in the case $p = 2$, particularly useful for applications. There are two methods for the numerical solutions of the BVPs for harmonic functions which, in fact, are founded on the completeness properties of the harmonic polynomials in the case $p = 2$}. We refer to 
\cite[pp.304--305]{fichera1979}
 for a description of these numerical methods.

  We observe an aspect of the results contained in the present paper. 
  Let us consider a BVP for which a certain completeness property holds on the boundary of a particular domain.
  In all the results obtained so far,  the fact that this completeness holds or not on other domains seems to depend exclusively on the topology. For example, in the case of the Laplace equation, harmonic polynomials are complete on the boundary if and only if the complement
  of the domain is connected. See also \cite{ C2019}, where the closure (in $L^p$ norm) of the linear space generated by the system of elastothermostatics polynomials
  on the boundary of a multiple connected domain is
  determined.
 
  The results contained in the present paper show that
  the completeness on the boundary may also depend on the geometry of the domain and not only on the topology. Indeed, for the third problem, we prove that the 
 closure of the linear space generated by the considered system has a different nature depending on whether the domain 
 is axially symmetric or not. 
 
 The paper is organized as follows. After some preliminaries given in Section  \ref{sec:prel}, 
 we introduce and study some properties of a particular class $\A^p(\Om)$ of solutions of elasticity system in Section \ref{sec:class}. 
 This  class  can be described as the class of solution of elasticity system for which displacements and stresses exist on the boundary in a weak sense and
 belong to $L^p$ spaces.  
 
 Section  \ref{sec:polyn} is devoted to recalling how to construct a system of elastic polynomials.
  We prove also a necessary and sufficient condition under which 
 two $L^p$ vectors defined on the boundary   can be considered as the displacements and the stresses of an element of $\A^p(\Om)$.
 
 In Section \ref{sec:Lp} we prove completeness theorems related to the III and IV BVPs in $L^p$ norm. The more
 delicate case of the completeness in uniform norm is treated in Section \ref{sec:C0}.

 \section{Preliminaries}\label{sec:prel}
 
 Let $\Om$ be a bounded domain in $\R^3$. If not otherwise specified, the domain $\Om$ is such that $\R^3\setminus \Om$ is connected 
 and its boundary $\Si$ is a $C^{2,\al}$ manifold ($0<\al\leq 1$). i.e. its curvature satisfies a H\"older condition
 of order $\al$ (see \cite[p.79]{kupra}).
 
 We denote by $L^p(\Om)\left(L^p(\Si)\right)$ the vector space of all measurable real functions $u$ such that $\left|u\right|^p$ is integrable over $\Om$ (over $\Si$) ($p \geq 1)$.
 
 If $k$ is any non-negative integer $C^k(\Om)\left(C^k(\overline{\Om})\right)$ indicates the class of all the  functions $u$ which are continuously differentiable up to the order $k$ in $\Om$ (in $\overline{\Om})$. 
 
 $\Cspt^{\infty}(\Om)$ is the space of  $C^{\infty}$ functions having compact support in $\Om$. 
 
$C^\beta(\overline{\Om})$ ($C^\beta(\Si)$) denotes the vector space of all continuous functions satisfying on $\overline{\Om}$ (on $\Si$) a uniform H\"older condition of some exponent $\beta$ ($0<\beta \leq1)$.

$C^{k,\beta}(\overline{\Om})$ stands for the sub-class of $C^k(\overline{\Om})$ consisting of functions $u$ such that $D^\alpha u \in C^\beta(\overline{\Om})$ ($|\alpha|=k$). The class $C^{1,\beta}(\Si)$ is defined in the same way, the derivatives with respect to the Cartesian coordinates being replaced by the tangential derivatives on $\Si$ (see \cite[p.80]{kupra}).

We shall consider also the Sobolev space on the boundary $W^{1,p}(\Si)$ (see, e.g, \cite[p.327]{kufner}).

 We denote by $E$ the operator of the classical  linear elasticity for a homogeneous isotropic body, i.e.
 $$
 Eu = \mu \De u +(\la+\mu) \nabla \dive u
 $$
 where $u$ is the vector $(u_1,u_2,u_3)$ and $\la$ and $\mu$ are the Lam\'e constants.  As usual (see, e.g., \cite[p.27]{kupra}), we assume that
 $$
 \mu>0\, , \quad  3\la+2\mu >0\, .
 $$
 
 By $Tu$ we denote the stress operator on the boundary of $\Om$:
 $$
 Tu = 2\mu \frac{\de u}{\de \nu} + \la(\dive u)\nu + \mu(\nu \wedge \curl u),
 $$
$\nu$ being the  outward unit normal to $\Si$.

Let  $\{\Ga_{ij}\}$ be the fundamental matrix:
$$
 \Ga_{ij}(x)=-\frac{\delta_{ij}}{4\pi\mu\, |x|} + \mu' \frac{\de^2 |x|}{\de x_i \de x_j}\, ,
 $$
 where
 $$
 \mu'=\frac{\la+\mu}{8\pi \mu (\la+2\mu)}\, .
 $$
We denote by $\Ga_i$  the vectors whose components are the elements of the $i$-th row of   $\{\Ga_{ij}\}$.

We recall that if $u$ and $v$ belong to $[C^{1}(\overline{\Om})\cap
 C^{2}(\Om)]^3$ and $E u$, $Ev$ are in $[L^{1}(\Om)]^3$, then
Betti's formula holds
 \begin{equation}\label{eq:betti}
\int_{\Om}  (u\cdot E v -   v \cdot E u)\, dx
=
\int_{\Si} (u\cdot Tv  - v \cdot Tu)\, d\si.
\end{equation}

We have also that if $u$ is a smooth vector in $[C^{1}(\overline{\Om})\cap
 C^{2}(\Om)]^3$ with $E u \in L^{1}(\Om)$,  the following ``Somigliana formula'' holds
 \begin{equation}\label{eq:somigliana}
\begin{gathered}
\int_{\Om} Eu(y)\cdot \Ga_i(x-y)\, dy + \int_{\Si}[u(y) \cdot T_y[\Ga_i(x-y)] -    \Ga_{i}(x-y)\cdot Tu(y)]\, d\si_{y} \\
=
\begin{cases}
u_i(x),  & x\in \Om, \\ 
 0  & x\in  \R^3 \setminus \overline{\Om}.
 \end{cases}
\end{gathered}
\end{equation}

We note that, as a consequence of this representation formula, we have that, if 
$w\in [\Cspt^{\infty}(A)]^3$,  $A$ being a domain in $\R^{3}$, then
\begin{equation}\label{eq:suppcomp}
w_i(x)=\int_{A} \Ga_{i}(x-y) \cdot Ew(y)\, dy, \quad \forall\ x\in \R^3.
\end{equation}

We recall the following jump formulas for the double layer potential
(see \cite[p.302]{kupra}) with density $\vf\in L^p(\Si)$ ($p> 1$)
\begin{equation}\label{eq:jumpdouble}
\lim_{x\to x_{0}^{\pm}} \int_{\Si} \vf(y) \cdot T_y[\Ga_i(x-y)] \, d\si_{y} =\pm \frac{1}{2}\vf_{i}(x_{0}) +  \int_{\Si} \vf(y) \cdot T_y[\Ga_i(x_{0}-y)] \, d\si_{y}
\end{equation}
for almost every $x_{0}\in\Si$. Here the $\lim_{x\to x_{0}^{+}}$ ($\lim_{x\to x_{0}^{-}}$) is understood as an 
internal (external) angular boundary value. We note also that the integral on the right-hand side is a singular integral. 
The presence of the factor $1/2$ and the differences in the sign with respect to the formula given in \cite{kupra}
are due to a different normalization in the definition of fundamental solution. 

We have also
\begin{equation}\label{eq:jumpdersimple}
\lim_{x\to x_{0}^{\pm}} \int_{\Si} \vf(y) \cdot T_{x_{0}}[\Ga_i(x-y)] \, d\si_{y} =\mp \frac{1}{2}\vf_{i}(x_{0}) +  \int_{\Si} \vf(y) \cdot T_{x_{0}}[\Ga_i(x_{0}-y)] \, d\si_{y}
\end{equation}
a.e. on $\Si$.

According to \cite{kupra}, the III boundary value problem is 
\begin{equation}\label{eq:BVPIII}
\begin{cases}
Eu=0\, , & \text{ in  }  \Om\, , \\ 
 u\cdot \nu = \vf\, , \ 
 \ Tu - (Tu\cdot \nu)\nu = \Phi\, ,& \text{ on  }  \Si 
\end{cases}
\end{equation}
and the IV is
\begin{equation}\label{eq:BVPIV}
\begin{cases}
Eu=0  & \text{ in  }  \Om\, , \\ 
 u- (u\cdot \nu)\nu = \Psi\, , \ 
 \ Tu\cdot \nu = \psi\, ,& \text{ on  }  \Si 
\end{cases}
\end{equation}

We remark that if $\Ga$ is a vector, then $\Ga - (\Ga\cdot \nu)\nu$ is orthogonal  to $\nu$ on $\Si$. This implies that
 if there exists a solution of the BVP \eqref{eq:BVPIII} or  of the BVP \eqref{eq:BVPIV} , we must have $\Phi\cdot\nu =0 $  and  $\Psi\cdot\nu =0 $ on $\Si$, respectively.
This is why we introduce the following  subspace of  $[C^{\be}(\Si)]^{3}$ 
 $$
 [C^{\be}(\Si)]^{3}_{0}=\{ F\in [C^\be(\Si)]^{3}\ |\ F \cdot \nu = 0 \ \text{ on } \Si\}
 $$
 and, analogously,
 $$
 [L^{p}(\Si)]^{3}_{0}=\{ F\in [L^{p}(\Si)]^{3}\ |\ F \cdot \nu = 0 \ \text{ a.e. on } \Si\}\, .
 $$
 We set also
 $$
 [C^{1,\be}(\Si)]^{3}_{0}= [C^{1,\be}(\Si)]^{3} \cap [C^{\be}(\Si)]^{3}_{0}\, ,
 \ \ 
 [W^{1,p}(\Si)]^{3}_{0}= [W^{1,p}(\Si)]^{3} \cap [L^{p}(\Si)]^{3}_{0}\, .
 $$
 We remark that in the problem \eqref{eq:BVPIII} we may have other compatibility conditions, according to the geometry
 of $\Om$. In fact, let $u$ be a regular solution of \eqref{eq:BVPIII}; by Betti's formula, we get
  $$
 \int_{\Si}Tu\cdot \om \, d\si = 0
 $$
 for any rigid displacement $\om$.  Moreover, since
 $$
 Tu\cdot \om = 
 (Tu - (Tu\cdot\nu)\nu)\cdot\om  + (Tu\cdot\nu)(\om\cdot\nu)
 $$
 we must have
 $$
 \int_{\Si} \Phi \cdot \om \, d\si = 0
 $$
 for any rigid displacement $\om$ which is tangential, i.e. such that 
 \begin{equation}\label{eq:omnu=0}
\om\cdot\nu=0 \quad \text{ on } \Si.
\end{equation}
If $\Si $ is not an axially symmetrical surface condition \eqref{eq:omnu=0} implies that the rigid displacement $\om$ vanishes. 
However, if $\Si$ is an axially symmetrical surface there exist non-vanishing tangential rigid displacements. 
In particular, if $\Si$ is a sphere, then
 a rigid displacement is tangential
 if and only if  
 $\om=b \wedge (x-x^0)$, where $b\in\R^3$ is arbitrary and $x^0$ is the center of the sphere.
  If $\Si$ is an axially symmetrical surface different from a sphere, then 
  a rigid displacement $\om$ satisfies \eqref{eq:omnu=0}
 if and only if  $\om = c_0\,  [a \wedge (x-x^0)]$, where  $c_0$ is an arbitrary real number, 
 $a$ is the unit vector of the rotation axis and $x^0$ is an arbitrary point on this axis.

  We have the following existence and uniqueness theorems for Problems III and IV.
 \begin{theorem}\label{th:exuniqIII}
    If $\Si$ is not an axially symmetrical surface, then 
    for any given $(\vf, \Phi)\in C^{1,\be}(\Si)\times [C^{\be}(\Si)]^{3}_{0}$ problem \eqref{eq:BVPIII} has a unique regular solution. If 
    $\Si$ is an axially symmetrical surface different from a sphere, then,   given  $(\vf, \Phi)\in C^{1,\be}(\Si)\times [C^{\be}(\Si)]^{3}_{0}$ such that
    \begin{equation*}
\int_{\Si} \Phi \cdot  [a \wedge (x-x^0)]\, d\si = 0
\end{equation*}
where $a$ is the unit vector of the rotation axis and $x^0$ is an arbitrary point on this axis, there exists a regular solution of problem \eqref{eq:BVPIII}.
Any two regular solutions may differ by an additive vector of the form
$$
c_0\,  [a \wedge (x-x^0)]\, ,
$$
$c_0$ being an arbitrary real number. Finally if $\Si$ is a sphere,  given  $(\vf, \Phi)\in C^{1,\be}(\Si)\times [C^{\be}(\Si)]^{3}_{0}$ such that
 \begin{equation*}
\int_{\Si} \Phi \cdot  [b \wedge (x-x^0)] \, d\si= 0, \quad \forall\ b\in \R^3,
\end{equation*}
$x^0$ being the center of the sphere, there exists a regular solution of problem \eqref{eq:BVPIII}. In this case
two regular solutions may differ by an additive vector of the form
$$
\sum_{k}c_k\, [a_k \wedge (x-x^0)] \, ,
$$
where $c_k\in \R$ and $a_k \in \R^3$ are arbitrary.
\end{theorem}

For the proof we refer to Theorems 1.8, p.115 and 5.12, p.378 of \cite{kupra}.

 \begin{theorem}\label{th:exuniqIV}
 For any given $(\Psi,\psi)\in  [C^{1,\be}(\Si)]^{3}_{0}\times C^{\be}(\Si)$ problem \eqref{eq:BVPIV} has a unique regular solution.
\end{theorem}

For the proof see  \cite[Th. 5.15, p.380]{kupra}.

\section{The class $\A^p(\Om)$}\label{sec:class}

In the present section, we introduce two classes of solutions of the Lam\'e system $Eu=0$. The first one, which we call $\A^p(\Om)$, is roughly speaking the class of the solutions $u$ of the Lam\'e system such that $u$ and $Tu$ do exist on the boundary in a weak sense and  belong to 
$[L^p(\Si)]^3$. More precisely
\begin{gather}
\A^p(\Om) = \{ u\in [L^1(\Om)]^3\ |\ \exists\ \al, \be \in [L^p(\Si)]^3:  \notag\\
\int_{\Om} u \cdot Ew\, dx = \int_{\Si}(\al\cdot Tw - \be \cdot w)\, d\si, \quad \forall\ w\in [C^{\infty}(\R^3)]^3\}.
\label{eq:defA}
\end{gather}
If not otherwise specified, we suppose $1< p<\infty$.

\begin{remark}\label{rem:1}
If $u\in \A^p(\Om)$, then $u\in [C^{\infty}(\Om)]^3$ and it satisfies $Eu=0$. This follows immediately from classical results,
since \eqref{eq:defA} implies
$$
\int_{\Om} u \cdot Ew\, dx = 0, \quad \forall\ w\in [\Cspt^{\infty}(\Om)]^3.
$$
\end{remark}

\begin{remark}\label{rem:2}
In condition \eqref{eq:defA} it is sufficient to consider $w$ which are
$C^{\infty}$ in a  neighborhood of $\overline{\Om}$.
\end{remark}

The other class, $\B^p(\Om)$, is defined as follows: we say that the vector $u$ belongs to  $\B^p(\Om)$ if there exists
$ \al, \be \in [L^p(\Si)]^3$ such that
\begin{equation}\label{eq:defB}
\int_{\Si}[\al(y) \cdot T_y[\Ga_i(x-y)] -    \be(y) \cdot \Ga_{i}(x-y) ]\, d\si_{y} =
\begin{cases}
u_i(x),  & x\in \Om, \\ 
 0  & x\in  \R^3 \setminus \overline{\Om}.
\end{cases}
\end{equation}

 Actually, these two classes coincide.
\begin{theorem}\label{th:=}
$\A^p(\Om)=\B^p(\Om)$.
 \end{theorem}
 \begin{proof}
 Suppose that $u$ belongs to $\A^p(\Om)$. If $x\notin \overline{\Om}$,
the vector  $\Ga_{i}(x -\cdot)$ is $C^{\infty}$ in a  neighborhood of $\overline{\Om}$ and  satisfies the system $Ew=0$.  Therefore,  \eqref{eq:defB} for $x\notin \overline{\Om}$ follows immediately from
\eqref{eq:defA} and Remark \ref{rem:2}. If $x$ is fixed in $\Om$, 
consider the scalar function $\vf_{\ep}(y)=\vf((x-y)/\ep)$, where $\vf \in C^{\infty}(\R^n)$
is such that $\vf(y)=0$ if $|y|\leq 1/2$ and $\vf(y)=1$ if $|y|\geq 1$.
Take the vector $w_{\ep}(y) = \vf_{\ep}(y)\, \Ga_{i}(x,y)$ 
where $\ep>0$ is such that $\overline{B_{\ep}(x)}\subset \Om$. 
Equation \eqref{eq:defA} leads to
\begin{gather*}
\int_{\Si}[\al(y) \cdot T_y[\Ga_i(x-y)] -    \Ga_{i}(x-y)\cdot \be(y)]\, d\si_{y} 
 =\int_{\Si}[\al \cdot Tw_{\ep} -    w_{\ep}\cdot \be]\, d\si\\
 = \int_{\Om} u \cdot Ew_{\ep}\, dx =  \int_{B_{\ep}(x)} u \cdot Ew_{\ep}\, dx\, .
\end{gather*}
The last equality holds because $Ew_{\ep}(y) = E_{y} [\Ga_{i}(x-y)]=0$ for $|x-y|>\ep$.
Keeping in mind \eqref{eq:betti}, \eqref{eq:somigliana} and Remark \ref{rem:1}, we can write the last integral as
\begin{gather*}
\int_{\de B_{\ep}(x)} (u\cdot Tw_{\ep} - w_{\ep} \cdot Tu )\, d\si
\\ =
\int_{\de B_{\ep}(x)} [u(y)\cdot T_y[\Ga_i(x-y)]  - \Ga_{i}(x-y) \cdot Tu(y) ]\, d\si_{y}=
u_{i}(x).
\end{gather*}
This completes the proof of \eqref{eq:defB}. 
 
 Conversely, if $u\in \B^p(\Om)$,
known results of potential theory show that $u$ belongs to $[L^{p}(\Om)]^3$.
Thanks to Remark \ref{rem:1}, it is enough to prove that \eqref{eq:defA} holds for any
$w\in [\Cspt^{\infty}(A)]^3$, where $A$ is a domain such that $\overline{\Om}\subset A$.
Equation \eqref{eq:defB} implies
\begin{gather*}
\int_{\Om}  u \cdot Ew\, dx =   
 \int_{\Om} \left( \int_{\Si}[\al(y) \cdot T_y[\Ga_i(x-y)] -    \Ga_{i}(x-y)\cdot \be(y)]\, d\si_{y} \right)\cdot Ew\, dx \\
 =  \int_{A} \left( \int_{\Si}[\al(y) \cdot T_y[\Ga_i(x-y)] -    \Ga_{i}(x-y)\cdot \be(y)]\, d\si_{y} \right)\cdot Ew\, dx \, .
\end{gather*}

Applying Fubini theorem, from \eqref{eq:suppcomp}
we get \eqref{eq:defA}.
 \end{proof}

 \begin{theorem}\label{th:abv}
    If $u\in \B^p(\Om)$ the following internal angular boundary values do exist for almost any $x_{0}\in \Si$
    and we have
    \begin{equation}\label{eq:abv}
\lim_{x\to x_0} u(x) = \al(x_{0})\, , \quad \lim_{x\to x_0} Tu(x) = \be(x_{0})\, .
\end{equation}
\end{theorem}
\begin{proof}
This result follows immediately from known jump formulas and Liapunov--Tauber theorem for elastic potentials (see \cite[Ch. V]{kupra}).
\end{proof}
 
 We note the following uniqueness theorems
 
 \begin{theorem}
    If $u\in \A^p(\Om)$  is such that $\al=0$ then $u= 0$ in $\Om$.
\end{theorem}
\begin{proof}
In view of equation \eqref{eq:defB} we can write
$$
-\int_{\Si} \be(y) \cdot \Ga_{i}(x-y)\, d\si_{y} =
\begin{cases}
u_i(x),  & x\in \Om, \\ 
 0  & x\in  \R^3 \setminus \overline{\Om}
 \end{cases}
$$
and \eqref{eq:jumpdersimple} leads to
$$
 \frac{1}{2}\be_{i}(x) +  \int_{\Si} \be(y) \cdot T_{x}[\Ga_i(x-y)] \, d\si_{y} =0
$$
for a.e. $x\in\Si$.  From known results on singular integral equations, we deduce that $\be$ is H\"older continuous.
This implies that $u$ is a regular solution of the Dirichlet problem and the thesis follows from a known uniqueness 
result (see \cite[p.115]{kupra}).
\end{proof}
 
 With the same proof, we can prove the uniqueness (mod. a rigid displacement) theorem for the second BVP.
  \begin{theorem}\label{th:neu!}
    If $u\in \A^p(\Om)$  is such that $\be=0$ then $u$ is a rigid displacement.
\end{theorem}
\begin{proof}
Because of equation \eqref{eq:defB} we can write
$$
\int_{\Si} \al(y)  \cdot T_y[\Ga_i(x-y)] \, d\si_{y} =
\begin{cases}
u_i(x),  & x\in \Om, \\ 
 0  & x\in  \R^3 \setminus \overline{\Om}
 \end{cases}
$$
and \eqref{eq:jumpdouble} leads to
$$
 -\frac{1}{2}\al_{i}(x) +  \int_{\Si} \al(y) \cdot T_{x}[\Ga_i(x-y)] \, d\si_{y} =0
$$
for a.e. $x\in\Si$.  From known results on singular integral equations, we deduce that $\al$ is H\"older continuous.
This implies that $u$ is a regular solution of the second BVP problem and the thesis follows from a known uniqueness 
result (see \cite[p.115]{kupra}).
\end{proof}

 \begin{theorem}\label{th:repres}
    The vector $u$ belongs to $\A^p(\Om)$ if and only if there exists $\vf\in [L^p(\Si)]^3$ such that
\begin{equation}\label{eq:slp}
u_i(x) = \int_{\Si} \vf(y) \cdot \Ga_{i}(x-y)\, d\si_y\, , \quad x\in\Om\, .
\end{equation}
\end{theorem}
\begin{proof}
Let $u$ be the simple layer potential \eqref{eq:slp} and let $w$ be an arbitrary vector in  $[\Cspt^{\infty}(\R^3)]^3$. We have
\begin{equation}\label{eq:conto1}
\int_{\Om} u \cdot Ew\, dx = \int_{\Si}\vf_{i}(y)\, d\si_{y} \int_{\Om} E_{j}w(x) \Ga_{ij}(x-y)\, dx \, .
\end{equation}

In view of Somigliana's formula \eqref{eq:somigliana} we get
$$
\int_{\Om} Ew(x)\cdot \Ga_i(y-x)\, dx= - \int_{\Si}[w(x) \cdot T_x[\Ga_i(y-x)] -    \Ga_{i}(y-x)\cdot Tw(x)]\, d\si_{x} +w_{i}(y)
$$
for any $y\in\Om$. Considering the internal angular boundary value for $y\to y_{0}$, $y_{0}$ being a point on $\Si$, 
and applying \eqref{eq:jumpdouble}, we obtain
\begin{gather*}
\int_{\Om} Ew(x)\cdot \Ga_i(y_{0}-x)\, dx \\
=\frac{1}{2} w_{i}(y_{0})- \int_{\Si}[w(x) \cdot T_x[\Ga_i(y_{0}-x)] -    \Ga_{i}(y_{0}-x)\cdot Tw(x)]\, d\si_{x}
\end{gather*}
a.e. on $\Si$.   Putting this equality in \eqref{eq:conto1} we find
\begin{gather*}
\int_{\Om} u \cdot Ew\, dx \\
= 
 \int_{\Si}\vf_{i}(y)\left( \frac{1}{2} w_{i}(y)- \int_{\Si}[w(x) \cdot T_x[\Ga_i(y-x)] -    \Ga_{i}(y-x)\cdot Tw(x)]\, d\si_{x}\right) 
 d\si_{y} \\
 =  \int_{\Si}  \biggl[w_{i}(x) \left( \frac{1}{2}\vf_{i}(x) - \int_{\Si} \vf (y)\cdot T_x[\Ga_i(y-x)]\, d\si_{y}\right)\\
 + T_{j}w(x)  \int_{\Si}\vf _{i}(y)\Ga_{ij}(y-x)\, d\si_{y}\biggr] =
\int_{\Si}(\al\cdot Tw - w\cdot \be)\, d\si,
 \end{gather*}
 where
 $$
 \al_{i}(x)=   \int_{\Si}\vf _{j}(y)\Ga_{ij}(y-x)\, d\si_{y}\, , \
 \be_{i}(x) =  -\frac{1}{2}\vf_{i}(x) + \int_{\Si} \vf (y)\cdot T_x[\Ga_i(x-y)]\, d\si_{y}\, .
 $$
 
The vectors $\al$ and $\be$ belonging to $[L^p(\Si)]^3$,  we have proved that $u\in\A^p(\Om)$.

Conversely, let $u\in \A^p(\Om)$. Consider the singular integral equation
\begin{equation}\label{eq:singint}
- \frac{1}{2}\vf_{i}(x) +  \int_{\Si} \vf(y) \cdot T_{x}[\Ga_i(x-y)] \, d\si_{y} = \be(x)\, .
\end{equation}

 From \eqref{eq:defA} it follows 
 $$
 \int_{\Si} \be \cdot (a + b\w x)\, d\si =0, \quad \forall\ a,b\in \R^3\, .
 $$
 This shows that there exists a solution $\vf$ of \eqref{eq:singint}.  The difference
 $$
 u_i(x) - \int_{\Si} \vf(y) \cdot \Ga_i(x-y) \, d\si_{y}
 $$
 belongs to $\A^p(\Om)$ and by  Theorem \ref{th:neu!} it is a rigid displacement. Since a rigid displacement can be written as a simple layer potential, the theorem is proved.
\end{proof}

 Natroshvili \cite{natro1988} proved that regular solutions of Problems III and IV given in Theorems 
 \ref{th:exuniqIII} and \ref{th:exuniqIV} can be represented by simple layer potentials. 
 He obtains this representation by writing a singular integral system for the density of such potentials
 which is equivalent to the relevant BVP. For the details, we refer to   \cite{natro1988}.
 It is straightforward to consider these singular integral systems in $L^p$  and to obtain similar results when the data belong to $L^p(\Si)$. Combining this fact with  Theorem \ref{th:repres},
 we have immediately the
 following results

\begin{theorem}\label{th:exuniqIIILp}
    If $\Si$ is not an axially symmetrical surface, then 
    for any given $(\vf, \Phi)\in W^{1,p}(\Si)\times [L^{p}(\Si)]^{3}_{0}$ problem \eqref{eq:BVPIII} has a unique  solution
    in $\A^{p}(\Om)$. If 
    $\Si$ is an axially symmetrical surface different from a sphere, then,   given  $(\vf, \Phi)\in W^{1,p}(\Si)\times [L^{p}(\Si)]^{3}_{0}$ such that
    \begin{equation*}
\int_{\Si} \Phi \cdot  [a \wedge (x-x^0)]\, d\si = 0
\end{equation*}
where $a$ is the unit vector of the rotation axis and $x^0$ is an arbitrary point on this axis, there exists a  solution of  problem \eqref{eq:BVPIII} in $\A^{p}(\Om)$.
Any two regular solutions may differ by  an additive vector of the form
$$
c_0\,  [a \wedge (x-x^0)]\, ,
$$
$c_0$ being an arbitrary real number. Finally if $\Si$ is a sphere,  given  $(\vf, \Phi)\in W^{1,p}(\Si)\times [L^{p}(\Si)]^{3}_{0}$ such that
 \begin{equation*}
\int_{\Si} \Phi \cdot  [b \wedge (x-x^0)] \, d\si= 0, \quad \forall\ b\in \R^3,
\end{equation*}
$x^0$ being the center of the sphere, there exists  a  solution of problem \eqref{eq:BVPIII} in $\A^{p}(\Om)$.
 In this case
two regular solutions may differ by an additive vector of the form
$$
\sum_{k}c_k\, [a_k \wedge (x-x^0)] \, ,
$$
where $c_k\in \R$ and $a_k \in \R^3$ are arbitrary. The solutions can be represented by simple layer potentials with densities in $ [L^{p}(\Si)]^{3}$.
\end{theorem}

 \begin{theorem}\label{th:exuniqIVLp}
 For any given $(\Psi,\psi)\in  [W^{1,p}(\Si)]^{3}_{0}\times [L^{p}(\Si)]$, problem \eqref{eq:BVPIV} has a unique solution
 in $\A^{p}(\Om)$, which can be represented by a simple layer potential with density in $[L^{p}(\Si)]^{3}$.
\end{theorem}

  \section{Polynomial solutions}\label{sec:polyn}
 
 In this section, we recall how to construct a complete system $\{p_{k}\}$ of polynomial solutions 
 of the Lam\'e system 
 \begin{equation}\label{eq:eqhom}
Eu =0\, .
\end{equation}
 This means that any polynomial solution of \eqref{eq:eqhom} - which we shall call
 elastic (vector) polynomials - 
can be written as a finite linear combination of elements of $\{p_{k}\}$.

 The following theorem is proved in \cite[p.59]{Fichera50}
 \begin{theorem}\label{th:hompol}
   A vector $u$ homogeneous of degree $k\in \Z$ is solution of  system \eqref{eq:eqhom} if and only if there exists
   a harmonic vector $w$, homogeneous of the same degree $k$, such that
   $$
   u = w - \frac{\la+\mu}{2}\, \frac{|x|^2\, \nabla(\dive w)}{\la(k-1)+\mu(3k-2)}\, .
   $$
\end{theorem}


Let us denote by $\{\omega_{ks}\} (s = 1,\ldots , 2k + 1; k = 0, 1, \ldots)$ a complete
system of harmonic polynomials, i.e.
\begin{equation}\label{eq:harmpol}
\omega_{ks}(x)= |x|^{k}Y_{ks}
\left(\frac{x}{|x|}\right) \qquad (s=1,\ldots , 2k + 1; k=0,1,\ldots)
\end{equation}
$\{Y_{ks}\}$
being the system of spherical harmonics.

\begin{theorem}
  Let $v=(v_1,v_2,v_3)$ be a homogeneous polynomials of degree $k\in \N$. 
  The vector $v$ satisfies  elasticity system \eqref{eq:eqhom}
  if, and only if,  it  is a linear combination of the following
$6k+3$ polynomials:
\begin{equation}\label{eq:poly}
\begin{gathered}
\left( \om_{ks}\! +\Lambda_k |x|^2 \de_{11}\om_{ks},
 \Lambda_k |x|^2\de_{21}\om_{ks}, \Lambda_k |x|^2\de_{31}\om_{ks}\right) ,\\
\left(\Lambda_k |x|^2\de_{12}\om_{ks},  \om_{ks} + 
\Lambda_k  |x|^2\de_{22}\om_{ks},\Lambda_k  |x|^2\de_{32}\om_{ks},
   \right) , \\
\left(\Lambda_k  |x|^2\de_{13}\om_{ks}, 
 \Lambda_k  |x|^2\de_{23}\om_{ks}, \om_{ks} +\Lambda_k  |x|^2\de_{33}\om_{ks},
 \right) ,
\end{gathered}
\end{equation}
$(s=1,\ldots,2k+1)$, where  
$\Lambda_k= -  \frac{\la+\mu}{2(\la(k-1)+\mu(3k-2))}$ and
$\omega_{ks}$ are the harmonic polynomials \eqref{eq:harmpol}.
\end{theorem}
\begin{proof}
The result follows immediately from Theorem \ref{th:hompol} and from the well-known fact that 
any harmonic polynomial can be written as a linear combination of 
elements of $\{\omega_{ks}\}$.
\end{proof}

We denote by $\{p_{k}\}$ $(k=0,1,2\ldots$) the system 
 constituted by the vector polynomials  \eqref{eq:poly} (s=1,\ldots , 2k + 1; k=0,1,\ldots), 
ordered in one sequence.
Clearly any polynomial solution of elasticity system \eqref{eq:eqhom}
can be written as a finite linear combination of elements of $\{p_{k}\}$.

\begin{theorem}\label{th:Appolyn}
    Given $\al, \be \in L^p(\Si)$ there exists $u\in\A^p(\Om)$ satisfying \eqref{eq:defA} 
if and only if
\begin{equation}\label{eq:condort}
 \int_{\Si}(\al\cdot Tp_{k} - \be \cdot p_{k})\, d\si =0\, , \quad k=0,1,2,\ldots
\end{equation}
\end{theorem}
\begin{proof}
The necessity of \eqref{eq:condort} follows immediately from \eqref{eq:defA}. Conversely, let us prove  first that
\begin{equation}\label{eq:atfirst}
\int_{\Si}[\al(y) \cdot T_y[\Ga_i(x-y)] -    \be(y) \cdot \Ga_{i}(x-y) ]\, d\si_{y}=0, \quad \forall\, x\in \R^3\setminus \overline{\Om}.
\end{equation}

It is known that the following expansions hold (for a proof, see \cite[p.79]{CiNino})
\begin{gather*}
\frac{1}{|x-y|}=\sum_{k=0}^{\infty}
\sum_{s=0}^{2k+1}\lambda_{ks}(x)\, \omega_{ks}(y), \\
|x-y|= \sum_{k=0}^{\infty}
\sum_{s=0}^{2k+1}
\left[\chi_{ks}(x)\, \omega_{ks}(y) + 
\tau_{ks}(x)\, |y|^{2}\omega_{ks}(y)\right]
\end{gather*}
where $\omega_{ks}$ are the harmonic polynomials \eqref{eq:harmpol}
and
\begin{gather*}
\lambda_{ks}(x)= \frac{1}{|x|^{k+1}}Y_{ks}
\left(\frac{x}{|x|}\right), \\
\chi_{ks}(x) = - \frac{1}{2k-1}\frac{1}{|x|^{k-1}}Y_{ks}
\left(\frac{x}{|x|}\right), \quad
\tau_{ks}(x) = \frac{1}{2k+3}\frac{1}{|x|^{k+1}}Y_{ks}
\left(\frac{x}{|x|}\right).
\end{gather*}
If we fix $x\neq 0$, these expansions and their derivatives uniformly converge in
any compact set contained in the ball $|y|<|x|$. 
Observing that $(2k+3)\tau_{ks}(x) =\lambda_{ks}(x)$, we may write
\begin{equation}\label{eq:gexp}
\begin{gathered}
\Ga_{ij}(x-y)
=  \mu' \sum_{k=0}^{\infty}
\sum_{s=0}^{2k+1}  \chi_{ks}(x) \frac{\partial^{2}}{\partial y_{i}\partial y_{j}}\omega_{ks}(y)\\
+ \sum_{k=0}^{\infty}
\sum_{s=0}^{2k+1} \lambda_{ks}(x)
\left(- \frac{\delta_{ij}}{4\pi\mu}\, \omega_{ks}(y)
+ \frac{\mu'}{2k+3}\, \frac{\partial^{2}}{\partial y_{i}\partial y_{j}}(|y|^{2}\omega_{ks}(y))\right).
\end{gathered}
\end{equation}

Since the vectors
$$
\nabla \frac{\partial}{\partial y_{i}}\omega_{ks}(y)
$$
are harmonic and divergence-free, they satisfy the Lam\'e system $Eu=0$. 

Let us fix now $k$ and $s$ and denote by $\Phi_{ij}(y)$ 
the function
$$
\Phi_{ij}(y) =  -\frac{\delta_{ij}}{4\pi\mu}\, \omega_{ks}(y)
+ \frac{\mu'}{2k+3}\, \frac{\partial^{2}}{\partial y_{i}\partial y_{j}}(|y|^{2}\omega_{ks}(y))
$$
and by $\Phi_i$ the vector $(\Phi_{i1},\Phi_{i2}, \Phi_{i3})$. Since 
$\Delta(|y|^{2}\omega_{ks}(y))=2(2k+3)\omega_{ks}(y)$, we find
$$
\Delta \Phi_i(y) = 2\mu' \nabla  \frac{\partial}{\partial y_{i}}\omega_{ks}(y)
$$
and
$$
\dive \Phi_i(y)= \left(- \frac{1}{4\pi\mu} +2 \mu'\right)  \frac{\partial}{\partial y_{i}}\omega_{ks}(y)=
-\frac{1}{4\pi(\la+2\mu)}\, \frac{\partial}{\partial y_{i}}\omega_{ks}(y)\, .
$$

Therefore
$$
\mu \De \Phi_i + (\la+\mu) \nabla(\dive \Phi_i) = 
\left( 2\mu \mu' - \frac{\la+\mu}{4\pi(\la+2\mu)}\right) \nabla \frac{\partial}{\partial y_{i}}\omega_{ks}\, .
$$

Since
$$
2\mu \mu' = \frac{\la+\mu}{4\pi(\la+2\mu)}
$$
we deduce
$$
\mu \De \Phi_i + (\la+\mu) \nabla(\dive \Phi_i) =0.
$$

Because of \eqref{eq:gexp},  we have proved that
\begin{equation}\label{eq:expans}
\Ga_{i}(x-y)=\sum_{h=0}^{\infty}c_{ih}(x)w_{h}(y)
\end{equation}
uniformly for $y \in \overline{\Omega}$, provided that $|x| > R$,
where $R > \max_{y\in\overline{\Omega}} |y|$. Here $c_{ih}$ are 
scalar
functions and $w_{h}$
are vector homogenous polynomials satisfying the Lam\'e system $Ew=0$.

Let us consider now the potentials
$$
u_{i}(x) = \int_{\Si}[\al(y) \cdot T_y[\Ga_i(x-y)] -    \be(y) \cdot \Ga_{i}(x-y) ]\, d\si_{y}\, .
$$

Let us fix $R>\max_{y\in\overline{\Omega}} |y|$. In view of \eqref{eq:expans}, we have
$$
u_i(x)=\sum_{h=0}^{\infty}  c_{ih}(x)  \int_{\Si}(\al\cdot Tw_{h} - \be \cdot w_{h})\, d\si 
$$
for any $x$ such that $|x|>R$.  The vector polynomials $w_h$ are homogeneous and satisfy the Lam\'e system $Ew=0$.
Therefore
they can be written as a finite linear combination of $p_k$ and 
in view of \eqref{eq:condort} we get
$$
u_i(x)=0, \quad \forall\ x\, :\, |x|>R.
$$

Because of the analyticity of the potentials $u_i$  in the connected domain $\R^3\setminus \overline{\Om}$, we find
$$
u_i(x)=0, \quad \forall\ x\in \R^3\setminus \overline{\Om}
$$
and \eqref{eq:atfirst} is proved. This means that $u\in \B^p(\Om)$ and the result follows from Theorem \ref{th:=}.
\end{proof}

 \section{Completeness Theorems in $L^p$ norm}\label{sec:Lp}
 
 Let us prove the completeness theorems for the Problem III (see \eqref{eq:BVPIII}).
 \begin{theorem}\label{th:complIII}
    Let $1\leq p<\infty$ and let $\Si$ be a not axially symmetrical surface. The system 
    \begin{equation}\label{eq:systemIII}
\left\{   p_k \cdot \nu \, , Tp_k - (Tp_k \cdot \nu)\nu
\right\}
\end{equation}
is complete in $L^p(\Si)\times [L^{p}(\Si)]_{0}^3$.
\end{theorem}
\begin{proof}
Suppose first $1<p<\infty$.
We have to show that if a functional 
$F\in (L^p(\Si)\times [L^{p}(\Si)]^3)^*$ vanishes on the manifold generated by
\eqref{eq:systemIII}, then it vanishes on $L^p(\Si)\times [L^{p}(\Si)]_{0}^3$. 

Let $(\psi,\Psi)\in L^q(\Si)\times [L^{q}(\Si)]^3$ ($1/p+1/q=1$) such that
\begin{equation}\label{eq:Psipsi}
\int_{\Si} [\psi\, p_k \cdot \nu + \Psi\cdot(Tp_k - (Tp_k \cdot \nu)\nu)]\, d\si = 0, \quad k=0,1,\ldots\, .
\end{equation}

This can be written as
$$
\int_{\Si} [(\Psi - (\Psi \cdot \nu)\nu)\cdot Tp_k + \psi\, \nu \cdot p_k ]\, d\si = 0, \quad k=0,1,\ldots\, .
$$

Thanks to Theorem \ref{th:Appolyn} the potential
\begin{equation}\label{eq:potentu}
\begin{gathered}
u_{i}(x) \\
= 
\int_{\Si} [(\Psi(y)- (\Psi(y) \cdot \nu(y))\nu(y))\cdot T_y[\Ga_{i}(x-y)] + \psi(y)\, \nu(y) \cdot \Ga_{i}(x-y)]\, d\si_{y}
\end{gathered} 
\end{equation}
belongs to $\A^{q}(\Om)$. Moreover (see \eqref{eq:abv}) we have
$$
u = \Psi - (\Psi \cdot \nu)\nu\, , \quad Tu= -\psi\, \nu
$$
a.e. on $\Si$, in the sense of the interior angular boundary values. This implies
\begin{equation}\label{eq:BVPIIIhom}
u\cdot \nu =0 , \quad Tu - (Tu \cdot \nu)\nu =0 
\end{equation}
a.e. on $\Si$. In other words the potential $u\in\A^q(\Om)$ is a solution of the homogeneous problem III. 
By Theorem \ref{th:exuniqIIILp} we have $u=0$ in $\Om$.  In view of Theorem \ref{th:abv} we get
$$
 \Psi - (\Psi \cdot \nu)\nu =0, \quad  \psi\, \nu =0
$$
a.e. on $\Si$. This implies
$$
\psi=0, \quad   \Psi = \ga\, \nu,
$$
where $\ga$ is a scalar function in $L^{q}(\Si)$.  Therefore
$$
\int_{\Si}(\psi\, f + \Psi\cdot F)\, d\si = \int_\Si \ga\, \nu\cdot F\, d\si =0
$$
for any $(f,F)\in L^p(\Si)\times [L^{p}(\Si)]_{0}^3$ and this completes the proof.

If $p=1$ one can repeat the proof, observing that since $(\psi,\Psi)$ is in  $L^{\infty}(\Si)\times [L^{\infty}(\Si)]^3$,
we have $(\psi,\Psi) \in L^s(\Si)\times [L^{s}(\Si)]^3$, for any $s>1$. Then we have just to replace  $\A^q(\Om)$ by $\A^s(\Om)$ for a fixed $s>1$. 
\end{proof}
 
If  $\Si$ is an axially symmetrical surface the system \eqref{eq:systemIII} is not complete in $L^p(\Si)\times [L^{p}(\Si)]_{0}^3$, but only in a certain subspace.
\begin{theorem}\label{th:14}
    Let $1\leq p<\infty$ and let $\Si$ be an axially symmetrical surface different from a sphere. Let  $a$ be the unit vector of the rotation axis and $x^0$ an arbitrary point on this axis. The system \eqref{eq:systemIII} is complete in the subspace $V^p$ of  
    $L^p(\Si)\times [L^{p}(\Si)]_{0}^3$ defined as follows
 $$
V^p=\left\{ (f,F)\in L^p(\Si)\times [L^{p}(\Si)]_{0}^3\ \left|\ \int_{\Si} F\cdot (a \wedge (x-x^0))\, d\si =0 \right. \right\}.
$$
\end{theorem}
\begin{proof}
Let $1<p<\infty$.
Let us remark that the system $\left\{   p_k \cdot \nu \, , Tp_k - (Tp_k \cdot \nu)\nu
\right\}$ is contained in $V^p$. Indeed, since $(a \wedge (x-x^0))\cdot \nu =0$ on $\Si$,  we have
$$
\int_{\Si}[Tp_k - (Tp_k \cdot \nu)\nu]\cdot  (a \wedge (x-x^0))\, d\si =
\int_{\Si}  Tp_k  \cdot  (a \wedge (x-x^0))\, d\si = 0
$$
by Betti's formula \eqref{eq:betti}. 

As in the proof or Theorem \ref{th:complIII},  we have that, if 
$(\psi,\Psi)\in L^q(\Si)\times [L^{q}(\Si)]^3$ ($1/p+1/q=1$) is such that conditions \eqref{eq:Psipsi} hold,
then the potential \eqref{eq:potentu} satisfies the homogeneous boundary conditions \eqref{eq:BVPIIIhom}. In the present
case, this implies that there exists a real constant $c$ such that
$$
u(x)= c\, (a \wedge (x-x^0)), \quad x\in \Om.
$$
This leads to
$$
 \Psi - (\Psi \cdot \nu)\nu =c\, (a \wedge (x-x^0)), \quad  \psi\, \nu =0
$$
a.e. on $\Si$, and then
$$
\psi=0, \quad   \Psi = \ga\, \nu + c\, (a \wedge (x-x^0)),
$$
a.e. on $\Si$,  $\ga$ being a scalar function in $L^{q}(\Si)$.  Therefore
$$
\int_{\Si}(\psi\, f + \Psi\cdot F)\, d\si = \int_\Si \ga\, \nu\cdot F\, d\si  +
c \int_{\Si}(a \wedge (x-x^0)) \cdot F\, d\si =0
$$
for any $(f,F)\in V^p$, which proves the theorem. 

The case $p=1$ can be treated as in Theorem \ref{th:complIII}.
\end{proof}
 
 Following the same lines we can prove the following result.
 \begin{theorem}
   Let $1\leq p<\infty$ and let  $\Si$ be a sphere centered at $x_0$. The system \eqref{eq:systemIII} is complete in the subspace $\widetilde{V}^p$ of  
    $L^p(\Si)\times [L^{p}(\Si)]_{0}^3$ defined as follows
 $$
\widetilde{V}^p=\left\{ (f,F)\in L^p(\Si)\times [L^{p}(\Si)]_{0}^3\ \left| \ \int_{\Si} F\cdot (b \wedge (x-x^0))\, d\si =0 \right.,\forall\ b\in\R^3 \right\}.
$$
\end{theorem}
\begin{proof}
The details are left to the reader.
\end{proof}

 Concerning the Problem IV (see \eqref{eq:BVPIV}) we have
  \begin{theorem}
    Let $1\leq p<\infty$. The system 
    \begin{equation}\label{eq:systemIV}
\left\{   p_k - (p_k \cdot \nu)\nu  \, , Tp_k \cdot \nu
\right\}
\end{equation}
is complete in $[L^{p}(\Si)]_{0}^3\times L^p(\Si)$.
\end{theorem}
 \begin{proof}
Let $1<p<\infty$. As in Theorem \ref{th:complIII}, we have to show that if a functional 
$F\in ([L^{p}(\Si)]^3 \times L^p(\Si))^*$ vanishes on the manifold generated by
\eqref{eq:systemIV}, then it vanishes on $[L^{p}(\Si)]_{0}^3\times L^p(\Si)$. 

Let $(\Psi,\psi)\in [L^{q}(\Si)]^3\times L^q(\Si)$ ($1/p+1/q=1$) such that
$$
\int_{\Si} [\Psi\cdot(p_k - (p_k \cdot \nu)\nu) + \psi\, Tp_k \cdot \nu]\, d\si = 0, \quad k=0,1,\ldots\, .
$$

This can be written as
$$
\int_{\Si} [\psi\, \nu \cdot Tp_k + (\Psi - (\Psi \cdot \nu)\nu)\cdot p_k ]\, d\si = 0, \quad k=0,1,\ldots\, .
$$

In view of Theorem \ref{th:Appolyn} the potential
\begin{gather*}
u_{i}(x) \\
= 
\int_{\Si} [\psi(y)\, \nu(y) \cdot T_y[\Ga_{i}(x-y)] + (\Psi(y)- (\Psi(y) \cdot \nu(y))\nu(y))\cdot \Ga_{i}(x-y) ]\, d\si_{y}
\end{gather*} 
belongs to $\A^{q}(\Om)$. This leads to 
$$
u =\psi\, \nu\, , \quad Tu= - \Psi + (\Psi \cdot \nu)\nu
$$
a.e. on $\Si$, and then
$$
u - (u \cdot \nu)\nu =0, \quad  Tu\cdot \nu =0 
$$
a.e. on $\Si$. In other words the potential $u\in\A^q(\Om)$ is a solution of the homogeneous problem IV. 
By Theorem \ref{th:exuniqIVLp} we have $u=0$ in $\Om$ and we deduce
$$
 \psi\, \nu=0 , \quad  \Psi - (\Psi \cdot \nu)\nu =0
$$
a.e. on $\Si$. This implies
$$
\psi=0, \quad   \Psi = \ga\, \nu,
$$
where $\ga$ is a scalar function in $L^{q}(\Si)$.  Therefore
$$
\int_{\Si}(\Psi\cdot F+ \psi\, f )\, d\si = \int_\Si \ga\, \nu\cdot F\, d\si =0
$$
for any $(F,f)\in  [L^{p}(\Si)]_{0}^3\times L^p(\Si)$.

The case $p=1$ can be treated as in Theorems \ref{th:complIII} and \ref{th:14}.
\end{proof}

 \section{Completeness in the uniform norm}\label{sec:C0}
 
 We first recall the following result
 proved in \cite[p.298--300]{C19882}
 
 \begin{theorem}\label{C1988}
    The system $\{T p_{k}\}$ is complete (in the uniform norm) in the space 
    of vectors $F\in [C^{0}(\Si)^{3}$ such that
    \begin{equation*}
\int_{\Si} F \cdot (a+b\wedge x)\, d\si = 0
\end{equation*}
for any $a,b \in \R^{3}$.
\end{theorem}

The next result shows that we can uniformly approximate  continuous
displacements and stresses simultaneously. Indeed let us consider
the subspace $\A^0(\Om) $ of $\A^p(\Om)$ given by the vectors $u$ such that $u$ and $Tu$ are continuous on $\Si$, i.e.
\begin{gather}
\A^0(\Om) = \{ u\in [L^1(\Om)]^3\ |\ \exists\ \al, \be \in [C^0(\Si)]^3:  \notag \\
\int_{\Om} u \cdot Ew\, dx = \int_{\Si}(\al\cdot Tw - \be \cdot w)\, d\si, \quad \forall\ w\in [C^{\infty}(\R^3)]^3\}.
\label{eq:defC0}
\end{gather}

It is easy to prove that  $\{(\al,\be)\ |\ u\in \A^0(\Om)\}$ equipped with the uniform norm is a Banach space.

\begin{theorem}\label{th:doppiacompl}
    The system $\{(p_{k},T p_{k})\}$ is complete (in the uniform norm) in the space $\{(\al,\be)\ |\ u\in \A^0(\Om)\}$.
\end{theorem}
\begin{proof}
Let  $(\al,\be)\in \{(\al,\be)\ |\ u\in \A^0(\Om)\}$; by \eqref{eq:defC0} we get
$$
\int_{\Si} \be \cdot w\, d\si = 0 
$$
for any rigid displacement $w$. Theorem \ref{C1988} shows that there exists a sequence of
elastic polynomials $\{\om_{n}\}$ such that $T\om_{n} \to \be$ in $[C^0(\Si)]^3$. On the other hand, for any fixed
number  $q>2$ there exists
a constant $K$ such that
$$
\inf _{a, b \in \R^3}\|\om+a+b \wedge x\|_{[C^0(\Om)]^3} \leq K\, \| T \om\|_{[L^q(\Si)]^3}
$$
for any sufficiently smooth $\om$ satisfying the system $E\om=0$ (see \cite[p.26]{C1990}).
This inequality shows that  $\{\om_{n}\}$ is a Cauchy sequence in the quotient of  $[C^0(\Si)]^3$ by the subspace spanned 
by rigid displacements. Since this subspace has finite dimension and the quotient space is complete, we can find a sequence of rigid displacements $\{\xi_n\}$ such that the sequence
$\{\widetilde{\om}_n=\om_n+\xi_n\}$ is uniformly convergent to a vector $\vf\in [C^0(\Si)]^3$.
We have
$$
\int_{\Si}(\widetilde{\om}_n \cdot Tp_k -  T\widetilde{\om}_n \cdot p_k)\, d\si = 0\,  , \qquad (k=0,1,\ldots)
$$
and then, letting $n\to \infty$,
$$
\int_{\Si}(\vf \cdot Tp_k -  \be \cdot p_k)\, d\si = 0\, , \qquad (k=0,1,\ldots) \, .
$$
Thanks to \eqref{eq:defC0} we can write
$$
\int_{\Si}(\vf -\al)\cdot Tp_k d\si = 0 \, , \qquad (k=0,1,\ldots) \, .
$$
Theorem \ref{C1988} implies that $\vf -\al$ coincides on $\Si$ with a rigid displacement $w_0$
and then  $\widetilde{\om}_n-w_0 \to \al$ uniformly on $\Si$. Since $T(\widetilde{\om}_n-w_0)=T \om_n$, we can say that
$(\widetilde{\om}_n-w_0, T(\widetilde{\om}_n-w_0)) \to (\al,\be)$ in the uniform norm and the proof is complete. 
\end{proof}

 We remark that, with a similar proof, we have also
 \begin{theorem}
    Let $1\leq p<\infty$. The system $\{(p_{k},T p_{k})\}$ is complete (in the $L^P$ norm) in the space $\{(\al,\be)\ |\ u\in \A^p(\Om)\}$.
\end{theorem}

 We can rephrase Theorem \ref{th:doppiacompl} in the following way, where we denote by $M(\Si)$ the dual space of $C^0(\Si)$, i.e. the space
 of real Borel measures defined on $\Si$.
\begin{theorem}\label{thmeasorth}
    If $(\mu,\ro)\in [M(\Si)]^3 \times [M(\Si)]^3$ is such that
    \footnote{If $\al=(\al_1,\al_2,\al_3)\in [C^0(\Si)]^3$ and $\mu=(\mu_1,\mu_2,\mu_3)\in  [M(\Si)]^3$, then
    $$
    \int_{\Si}  \al\, d\mu = \int_{\Si} \al_j\, d\mu_j\, .
    $$
    }
    \begin{equation*}
\int_{\Si} (p_k \, d\ro + Tp_{k}\, d\mu ) =0\, , \qquad (k=0,1,\ldots) \, ,
\end{equation*}
then 
\begin{equation*}
\int_{\Si} ( \al\, d\ro + \be \, d\mu) =0
\end{equation*}
for any $(\al,\be)\in \{(\al,\be)\ |\ u\in \A^0(\Om)\}$.
\end{theorem}
\begin{proof}
By the Hahn-Banach theorem  the statement is equivalent to the completeness proved in Theorem \ref{th:doppiacompl}.\end{proof}

 We are now in a position to prove the completeness theorems in the uniform norm for problems III and IV.
 
 \begin{theorem}\label{th:complIIIunif}
    Let $\Si$ be a not axially symmetrical surface. The system  \eqref{eq:systemIII}
is complete in $C^0(\Si)\times [C^{0}(\Si)]_{0}^3$\, .
\end{theorem}
\begin{proof}
We have to show that, if $(\ro_0, \mu)\in  M(\Si) \times [M(\Si)]^3$ is such that
\begin{equation}\label{eq:misIII1}
\int_{\Si} [p_k \cdot \nu \, d\ro_0 + (Tp_k - (Tp_k \cdot \nu)\nu)\, d\mu] = 0, \quad k=0,1,\ldots\, ,
\end{equation}
then $(\ro_0, \mu)$ vanishes on $C^0(\Si)\times [C^{0}(\Si)]_{0}^3$\, . This means that $\ro_0=0$ and
$$
\int_{\Si} F\, d\mu =0
$$
 for any $F\in [C^{0}(\Si)]_{0}^3$. We can rewrite \eqref{eq:misIII1} as
 $$
 \int_{\Si} [p_k \cdot \nu \, d\ro_0 + T_{j}p_{k}(\delta_{jh}-  \nu_j\nu_h)\, d\mu_{h}]=0
 \quad k=0,1,\ldots\, ,
 $$
 where $\mu=(\mu_1,\mu_2,\mu_3)$. By Theorem \ref{thmeasorth} we find
\begin{equation}\label{eq:misIII3}
\int_{\Si}[ \al \cdot \nu\, d\ro_{0} + \be_{j}(\delta_{jh}-  \nu_j\nu_h)\, d\mu_{h}]=0
\end{equation}
for any $(\al,\be)\in \{(\al,\be)\ |\ u\in \A^0(\Om)\}$.
Let now $(\vf,\Phi) \in C^{1,\la}(\Si) \times [C^{\la}(\Si)]^{3}_{0}$. Because of Theorem \ref{th:exuniqIII}, there exists a smooth solution $u$ of the Problem
\eqref{eq:BVPIII}. We have
\begin{gather*}
\int_{\Si} (\vf\, d\ro_0 + \Phi \, d\mu)= \int_{\Si} [u \cdot \nu \, d\ro_0 + (Tu - (Tu \cdot \nu)\nu)\, d\mu]\\
=\int_{\Si}[ u \cdot \nu\, d\ro_{0} + T_{j}u\, (\delta_{jh}-  \nu_j\nu_h)\, d\mu_{h}]\, .
\end{gather*}
Thanks to \eqref{eq:misIII3}, the last integral vanishes,  $(u,Tu)$ belonging to $\{(\al,\be)\ |\ u\in \A^0(\Om)\}$.
 We have proved that
\begin{equation}\label{eq:35}
\int_{\Si} (\vf\, d\ro_0 + \Phi \, d\mu)=  0
\end{equation}
for any $(\vf,\Phi) \in C^{1,\la}(\Si) \times [C^{\la}(\Si)]^{3}_{0}$. In particular, we have
$$
\int_{\Si} \vf\, d\ro_0=  0
$$
 for any $\vf\in C^{1,\la}(\Si)$ and then $\ro_0=0$.  Let $\Theta\in [C^{\la}(\Si)]^{3}$ and set
$\Phi=\Theta-(\Theta\cdot \nu)\nu$. Since $\Phi$ belongs to $[C^{\la}(\Si)]^{3}_{0}$, we get
 $$
0=\int_{\Si}\Phi\, d\mu = \int_{\Si}[\Theta-(\Theta\cdot \nu)\nu]\, d\mu =
\int_{\Si}\Theta_{j}(\delta_{jh}-  \nu_j\nu_h)\, d\mu_{h}\, .
$$
This implies
 $(\delta_{jh}-  \nu_j\nu_h)\, \mu_{h}=0$, i.e.
$\mu_{j}=\nu_j\nu_h\, \mu_{h}$.  Then
$$
\int_{\Si} F\, d\mu = \int_{\Si} F_j \,\nu_j\,\nu_h\, d\mu_{h} = 0
$$
for any $F\in [C^{0}(\Si)]_{0}^3$.
 
 \end{proof}
 
 \begin{theorem}\label{th:21}
    Let $\Si$ be an axially symmetrical surface different from a sphere. Let  $a$ be the unit vector of the rotation axis and $x^0$ an arbitrary point on this axis. The system \eqref{eq:systemIII} is complete in the subspace $V^0$ of  
    $C^0(\Si)\times [C^{0}(\Si)]_{0}^3$ defined as follows
 $$
V^0=\left\{ (f,F)\in C^0(\Si)\times [C^{0}(\Si)]_{0}^3\ \left|\ \int_{\Si} F\cdot (a \wedge (x-x^0))\, d\si \right. =0 \right\}.
$$
\end{theorem}
\begin{proof}
The argument used at the beginning of the proof of theorem \ref{th:14} shows that the system $\left\{ p_k \cdot \nu \, , Tp_k - (Tp_k \cdot \nu)\nu
\right\}$ is contained also in $V^0$. We have to show that, if $(\ro_0, \mu)\in  M(\Si) \times [M(\Si)]^3$ is such that
\eqref{eq:misIII1} holds, then $(\ro_0, \mu)$ vanishes on all $(f,F)\in C^0(\Si)\times [C^{0}(\Si)]_{0}^3$ such that
the orthogonality condition
 \begin{equation}\label{eq:orth1bis}
\int_{\Si} F \cdot  [a \wedge (x-x^0)]\, d\si = 0
\end{equation}
 holds.

As in Theorem \ref{th:complIIIunif}, we find that \eqref{eq:misIII3} holds  for any $(\al,\be)\in \{(\al,\be)\ |\ u\in \A^0(\Om)\}$.
Let us take $(\vf,\Phi) \in C^{1,\la}(\Si) \times [C^{\la}(\Si)]^{3}_{0}$ such that 
\begin{equation}\label{eq:orthgamma}
\int_{\Si} \Phi \cdot \gamma\, d\si = 0\, ,
\end{equation}
where $\gamma(x)= (a \wedge (x-x^0))/ \| a \wedge (x-x^0) \|_{[L^{2}(\Si)]^3}$. Let  $u$ be the smooth solution of the Problem
\eqref{eq:BVPIII} (see Theorem \ref{th:exuniqIII}). The same arguments used in the proof of Theorem \ref{th:complIIIunif} show that
\eqref{eq:35}
holds
for any $(\vf,\Phi) \in C^{1,\la}(\Si) \times [C^{\la}(\Si)]^{3}_{0}$ satisfying \eqref{eq:orthgamma} and in the same way we get $\ro_0=0$.
 Let $\Theta\in [C^{\la}(\Si)]^{3}$ and set
$\Phi=\Theta-(\Theta\cdot \nu)\nu-(\int_{\Si}\Theta\cdot \gamma d\si)\gamma$.  Since $\nu\cdot \gamma=0$, we have 
that $\Phi$ belongs to $[C^{\la}(\Si)]^{3}_{0}$ and satisfies \eqref{eq:orthgamma}.
Therefore
\begin{gather*}
0=\int_{\Si}\Phi\, d\mu = \int_{\Si}[\Theta-(\Theta\cdot \nu)\nu]\, d\mu -\int_{\Si} \Theta\cdot\gamma\, d\si \int_{\Si}\gamma\, d\mu\\
=
\int_{\Si}\Theta_{j}[(\delta_{jh}-  \nu_j\nu_h)\, d\mu_{h} - c\, \gamma_j\, d\si]\, ,
\end{gather*}
where
$$
c=\int_{\Si} \gamma\, d\mu\, .
$$

For the arbitrariness of $\Theta$ we find $(\delta_{jh}-  \nu_j\nu_h)\, \mu_{h} - c\, \gamma_j\, \si =0$, i.e.
$\mu_{j} = \nu_j\nu_h\, \mu_{h} +  c\, \gamma_j\, \si$. Then
$$
\int_{\Si} F\, d\mu = \int_{\Si} F_j \,\nu_j\,\nu_h\, d\mu_{h} + c \int_{\Si} F\cdot \gamma\, d\si= 0
$$
for any $F\in [C^{0}(\Si)]_{0}^3$ satisfying \eqref{eq:orth1bis} and the proof is complete.
 \end{proof}

 \begin{theorem}
    Let  $\Si$ be a sphere centered at $x_0$. The system \eqref{eq:systemIII} is complete in the subspace $\widetilde{V}^0$ of  
    $C^0(\Si)\times [C^{0}(\Si)]_{0}^3$ defined as follows
 $$
\widetilde{V}^0=\left\{ (f,F)\in C^0(\Si)\times [C^{0}(\Si)]_{0}^3\ \left| \ \int_{\Si} F\cdot (b \wedge (x-x^0))\, d\si =0 \right.,\forall\ b\in\R^3 \right\}.
$$
\end{theorem}
\begin{proof}
We only sketch the proof, which follows the lines of that of Theorem \ref{th:21}.  Let  $\gamma_j$ ($j=1,2,3$) be an orthonormal 
base (with respect the scalar product in $[L^{2}(\Si)]^3$) of the space $\{b \wedge (x-x^0) \, |\,  b\in\R^3\}$.
 Replace $V^0$ by $\widetilde{V}^0$ and 
in the last part of the proof, given $\Theta\in [C^{\la}(\Si)]^{3}$, define
$\Phi=\Theta-(\Theta\cdot \nu)\nu-(\int_{\Si}\Theta\cdot \gamma_j d\si)\gamma_j$. We omit the details.
\end{proof}

 As far as the completeness in the uniform norm for the fourth problem is concerned, we have the following result
 
  \begin{theorem}
   The system \eqref{eq:systemIV} is complete in $[C^{0}(\Si)]_{0}^3\times C^0(\Si)$.
\end{theorem}
 \begin{proof}
We have to show that, if $(\ro, \mu_0)\in [M(\Si)]^3\times M(\Si)$ is such that
\begin{equation}\label{eq:misure1}
\int_{\Si} [(p_k - (p_k \cdot \nu)\nu)\, d\ro + Tp_k \cdot \nu\, d\mu_0] = 0, \quad k=0,1,\ldots\, ,
\end{equation}
then $\mu_0=0$ and
$$
\int_{\Si} F\, d\ro =0
$$
 for any $F\in [C^{0}(\Si)]_{0}^3$. We can rewrite \eqref{eq:misure1} as
 \begin{equation*}
\int_{\Si} [ (p_{k,j}(\delta_{jh}-  \nu_j\nu_h)\, d\ro_{h}  + Tp_k \cdot \nu\, d\mu_0] = 0, \quad k=0,1,\ldots\, .
\end{equation*}
where $p_k=(p_{k,1},p_{k,2},p_{k,3})$ and $\ro=(\ro_1,\ro_2,\ro_3)$. By Theorem \ref{thmeasorth} we find
\begin{equation}\label{eq:misure3}
\int_{\Si}[ \al_{j}(\delta_{jh}-  \nu_j\nu_h)\, d\ro_{h} + \be \cdot \nu\, d\mu_0] =0
\end{equation}
for any $(\al,\be)\in \{(\al,\be)\ |\ u\in \A^0(\Om)\}$.

Let now $(\Psi,\psi) \in [C^{1,\la}(\Si)]^{3}_{0} \times C^{\la}(\Si)$ and let $u$ be the smooth solution of the Problem
\eqref{eq:BVPIV}. We have
\begin{gather*}
\int_{\Si} (\Psi\, d\ro + \psi\, d\mu_0) = 
\int_{\Si}[(u- (u\cdot \nu)\nu)\, d\ro +  Tu\cdot \nu\, d\mu_0] \\
= \int_{\Si}[ u_{j}(\delta_{jh}-  \nu_j\nu_h)\, d\ro_{h} + Tu\cdot \nu\, d\mu_0 ] \, .
\end{gather*}
The last integral vanishes in view of \eqref{eq:misure3},  $(u,Tu)$ belonging to $\{(\al,\be)\ |\ u\in \A^0(\Om)\}$.
We have proved that
$$
\int_{\Si} (\Psi\, d\ro + \psi\, d\mu_0) = 0
$$
for any $(\Psi,\psi) \in [C^{1,\la}(\Si)]^{3}_{0} \times C^{\la}(\Si)$. In particular, we have
$$
\int_{\Si} \psi\, d\mu_0 = 0
$$
for any $\psi \in C^{\la}(\Si)$, which implies $\mu_0=0$.  Let $\Theta\in [C^{1,\la}(\Si)]^{3}$ and set
$\Psi=\Theta-(\Theta\cdot \nu)\nu$. Since $\Psi$ belongs to $[C^{1,\la}(\Si)]^{3}_{0}$, we get
$$
0=\int_{\Si}\Psi\, d\ro = \int_{\Si}[\Theta-(\Theta\cdot \nu)\nu]\, d\ro =
\int_{\Si}\Theta_{j}(\delta_{jh}-  \nu_j\nu_h)\, d\ro_{h}\, .
$$
Because of the arbitrariness of $\Theta$ we deduce $(\delta_{jh}-  \nu_j\nu_h)\, \ro_{h}=0$, i.e.
$\ro_{j}=\nu_j\nu_h \ro_{h}$.  Then
$$
\int_{\Si} F\, d\ro = \int_{\Si} F_j \,\nu_j\,\nu_h\, d\ro_{h} = 0
$$
for any $F\in [C^{0}(\Si)]_{0}^3$ and the theorem is proved.
 \end{proof}

 \textbf{Acknowledgement.}
 A. Cialdea is member of Gruppo Nazionale per l’Analisi
Matematica, la Probabilit\`a e le loro Applicazioni (GNAMPA) 
of the Istituto Nazionale di Alta Matematica (INdAM)
and acknowledges the support from the project ``Perturbation problems and
asymptotics for elliptic differential equations: variational and potential
theoretic methods'' funded by the European Union - Next Generation EU and
by MUR ``Progetti di Ricerca di Rilevante Interesse Nazionale'' (PRIN) Bando
2022 grant 2022SENJZ3.


\end{document}